\documentclass[11pt]{amsart}
\usepackage{amssymb, latexsym}

 \theoremstyle{plain}
 \newtheorem{thm}{Theorem}[section]
 
 \newtheorem{lem}[thm]{Lemma}
 \newtheorem{prop}[thm]{Proposition}
 \newtheorem{ex}[thm]{Example}
 
 \newtheorem{lem-def}[thm]{Lemma/Definition}
 \newtheorem*{thm-no-num}{Theorem}
 \newtheorem*{cor-no-num}{Corollary}

 \newtheorem*{temp1}{\theoremname}
 \newcommand{\theoremname}{temp2}
 \newenvironment{no-num}[1]{\renewcommand{\theoremname}{#1}
    \begin{temp1}}
    {\end{temp1}}

 \theoremstyle{definition}
 \newtheorem{defn}[thm]{Definition}

 \theoremstyle{remark}
 \newtheorem{rem}[thm]{Remark}

 \numberwithin{equation}{section}

 \newcommand{\To}{\rightarrow}
 \newcommand{\bP}{\mathbb{P}}
 \newcommand{\bC}{\mathbb{C}}
 \newcommand{\bZ}{\mathbb{Z}}
 
 \newcommand{\bQ}{\mathbb{Q}}
 \newcommand{\bF}{\mathbb{F}}
 \newcommand{\bG}{\mathbb{G}}
 \newcommand{\bA}{\mathbb{A}}
 \newcommand{\bmu}{\mathbf{\mu}}

 \newcommand{\DM}{Deligne-Mumford}
 
 \newcommand{\sY}{\mathcal{Y}}
 \newcommand{\sX}{\mathcal{X}}
 \newcommand{\sC}{\mathcal{C}}
 \newcommand{\sK}{\mathcal{K}}
 \newcommand{\sM}{\mathcal{M}}
 
 \newcommand{\sU}{\mathcal{U}}
 
 \newcommand{\sO}{\mathcal{O}}
 \newcommand{\ii}{\mathcal{I}}

 \newcommand{\coloneq}{\mathrel{\mathop:}=}

 \newcommand{\ov}{\overline}

\usepackage{stmaryrd}
\usepackage{diagrams}
\newarrow{Eq}=====
\newarrow{Into}C--->
\newarrow{Mapsto}|--->

\begin{document}

\title[Frobenius Action on $\ell$-adic Chen-Ruan Cohomology]
 {Frobenius Action on $\ell$-adic Chen-Ruan Cohomology}
\author{ Michael A. Rose }
\address{Department of Mathematics, University of British Columbia, Vancouver, British Columbia, V6T1Z2,
 Canada}
\email{mrose@math.ubc.ca}
\date{September 13, 2007}

\begin{abstract}
We extend the definition of the Chen-Ruan cohomology ring to smooth, proper, tame, \DM\ stacks over fields of positive characteristic and prove that a modified version of the Frobenius action preserves the product.
\end{abstract}

\maketitle
\section{Introduction}
The Weil conjectures describe a strong relationship between the arithmetic and topological properties of a smooth projective scheme $X$ over Spec($\bZ$).  In \cite{Weil-1956}, Weil himself observed that the conjectures would follow from an appropriate cohomology theory for abstract schemes (analogous to singular cohomology for complex varieties), and in the 1960's a great amount of work was done by Artin, Deligne, Grothendieck and others to develop \emph{$\ell$-adic cohomology} for this purpose.  In particular, by considering a smooth reduction $X_{\bF_q}$ and applying the Lefschetz Trace Theorem to the geometric Frobenius morphism $F$ on $\ov{X} \coloneq X \times_{\bF_q} \ov{\bF}_q$ one obtains the fundamental equation
\begin{equation} \label{intro_eqn_1}
\text{det}(1 - F^*t\ |\ H^*(\ov{X}_{\acute{e}t}, \bQ_l))\ =\ \text{exp}\ \sum_{r = 1}^{\infty}\ |X(\bF_{q^r})| \frac{t^r}{r}\
\end{equation}
where \emph{det} denotes a graded determinant, and $l$ is coprime to $q$.

Now we replace $X$ with a smooth \DM\ stack $\sX$ over Spec($\bZ$) (further hypothesis to be considered below), and one is naturally led to consider the Weil Conjectures on $\sX$.  Again, one focuses on finding an appropriate cohomology theory.  On one hand, there is already a natural notion of $\ell$-adic cohomology in this setting, and several of the properties crucial to proving an analog of the Weil conjectures for stacks have been established: Poincar\'{e} Duality appeared in \cite{LO1, LO2} while a Lefschetz Trace Theorem was established in \cite{Be-Lef}.  On the other hand, while motivated by string theory, Chen and Ruan \cite{CR2} constructed a ring now bearing their names:  the \emph{Chen-Ruan cohomology ring}, $H^*_{CR}(\sX_{\bC})$.  (The authors actually defined the ring for an almost complex orbifold; the theory was developed for \DM\ stacks in \cite{AGV1, AGV2, AV1}).

The goal of this article is to study the relationship between the Chen-Ruan cohomology ring and the arithmetic properties of $\sX$.  If $\sX$ is now a smooth, tame, \DM\ $\bF_q$-stack with projective coarse moduli scheme and $\ov{\sX} \coloneq \sX \times_{\bF_q} \ov{\bF}_q$, then we define the \emph{$\ell$-adic Chen-Ruan cohomology ring} of $\ov{\sX}$ denoted $H^*_{CR}(\ov{\sX}_{\acute{e}t}, \bQ_l)$.  Most of the technical requirements for this construction already appeared in \cite{AGV2}.  Furthermore we construct an action of the arithmetic Frobenius on $H^*_{CR}(\ov{\sX}_{\acute{e}t}, \bQ_l)$.  (We use the arithmetic Frobenius as opposed to the geometric Frobenius merely to simplify the proofs.  In the context of Artin stacks however, this distinction is crucial for convergence issues.  See \cite{Be-Lef}).

This latter result deserves comment.  In general, $H^*_{CR}$ \emph{is not functorial}.  Let $F:\ov{\sX} \To \ov{\sX}$ denote the arithmetic Frobenius on $\ov{\sX}$.  While $F$ naturally induces a linear map $\ii_{\bmu}(F)^*: H^*_{CR}(\ov{\sX}_{\acute{e}t}, \bQ_l) \To H^*_{CR}(\ov{\sX}_{\acute{e}t}, \bQ_l)$ (see Section \ref{Arithmetic_Frobenius} for a detailed definition), $\ii_{\bmu}(F)^*$ is not a ring homomorphism.  However, the main proposition of this article shows that a slight modification of $\ii_{\bmu}(F)^*$ indeed preserves the product structure:
\begin{prop} \label{main}
The \emph{orbifold Frobenius morphism} $F_{orb}$ given by
\begin{alignat}{1}
F_{orb}: \ H^*_{CR}(\ov{\sX}_{\acute{e}t}, \bQ_l) &\longrightarrow H^*_{CR}(\ov{\sX}_{\acute{e}t}, \bQ_l) \notag \\
                                        \alpha       &\longmapsto q^{-\text{age}(\alpha)} \cdot \ii_{\bmu}(F_{\sX})^*(\alpha) \notag
\end{alignat}
is an isomorphism of graded rings.
\end{prop}
Here, $\text{age}: H^*_{CR}(\ov{\sX}_{\acute{e}t}, \bQ_l) \To \bQ$ is a function appearing also in the grading on $H^*_{CR}(\ov{\sX}_{\acute{e}t}, \bQ_l)$ described in Section \ref{l-adic_cohomology}.

We study the arithmetic information contained by this Galois representation in Section \ref{zeta} below.  In particular, we include the analogues of equation (\ref{intro_eqn_1}) with $H^*$ (resp. $F^*$) replaced by $H^*_{CR}$ (resp. $F^*_{CR}$), and we list a consequence of Yasuda's proof \cite{Yas} of an additive version of the Crepant Resolution Conjecture.

\emph{Conventions.}  Unless specified otherwise, assume all sheaves on a \DM\ stack $\sX$ are defined on the \'{e}tale site of $\sX$.  We fix an isomorphism of the Tate twist $\bQ_l(1) \cong \bQ_l$ as sheaves on $\ov{\sX}$, inducing an isomorphism $\bQ_l(r) \cong \bQ_l$ for each $r$.  Then for smooth $\sX$ of dimension $n$, the duality theorem in \cite[Theorem 7.7]{LO2} yields $H^i(\ov{\sX}, \bQ_l) \cong H^{2n-i}(\ov{\sX}, \bQ_l)^{\vee}$.  This isomorphism is used implicitly throughout.  We define homology groups $H_i(\ov{\sX}, \bQ_l) \coloneq H^{2n-i}(\ov{\sX}, \bQ_l)$.  To improve the exposition proofs of several lemmas appear only in the appendix.  Also, no attempt is made to make the statements as general as possible, and we work over explicit base fields $\bF_q$ and $\ov{\bF}_q$.

The author would like to thank Lev Borisov for inspiration and for encouraging me to pursue this subject, Andrei C\u{a}ld\u{a}rau and Jordan Ellenberg for useful conversations, and Hsian-Hua Tseng and Tom Graber for pointing out the relevance of Yasuda's work in this context.

\section{Arithmetic Frobenius and Inertia Stacks} \label{Arithmetic_Frobenius}
In this section, we review inertia stacks and the induced action of the arithmetic Frobenius on them.  The notion of inertia stack is essential to the Gromov-Witten theory of stacks, and it's $\ell$-adic cohomology will form the underlying vector space for main object of study in this article: the $\ell$-adic Chen-Ruan cohomology ring.

Let $\sX$ be a \DM\ stack over $\bF_q$.  Fix an algebraic closure $\ov{\bF}_q \supset \bF_q$ and denote $\ov{\sX} \coloneq \sX \times_{{\bF}_q} \ov{\bF}_q$.  Let $\ov{\bF}_q \xrightarrow{\phi} \ov{\bF}_q$ denote the Frobenius morphism given by $\lambda \To \lambda^q$, and let $F_{\ov{\sX}, q} \coloneq 1_{\sX} \times \ \text{Spec}(\phi)$ \ denote the \emph{arithmetic Frobenius} morphism on $\sX$.  The subscripts will be dropped when no confusion arises.

Let $\bmu_r \coloneq \text{Spec}(\bZ[t]/ \langle t^r - 1 \rangle)$ denote the group scheme over $\bZ$ of $\text{r}^{th}$ roots of unity, and also denote by $\bmu_r$ the base change to $\bF_q$ when the context is clear.  Let $B\bmu_r \coloneq [\text{Spec}(\bF_q) / \bmu_r]$ denote the quotient stack corresponding to the trivial action of $\bmu_r$ on $\text{Spec}(\bF_q)$.

\begin{defn} \label{inertia_defn}
  \hfill
  \begin{enumerate}
    \item We denote by $\ii_{\bmu_r}(\sX) \equiv \text{\underline{Hom}}^{\text{rep}}_{\ \bF_q}(B\bmu_r, \sX)$ the stack of representable 1-morphisms over $\bF_q$ from $B\bmu_r$ to $\sX$.\\
    \item The \emph{cyclotomic inertia stack} is given by
            $$\ii_{\bmu}(\sX) \coloneq \bigsqcup_r \ii_{\bmu_r}(\sX).$$
  \end{enumerate}
\end{defn}
$\ii_{\bmu}(\sX)$ is a \DM\ stack and is smooth when $\sX$ is tame and smooth (this follows from \cite[Section 3]{AGV2}).

\begin{ex} \label{inertia_example}
Let $G \To \text{Spec }(\bF_q)$ be a finite \'{e}tale group scheme acting on $\bF_q$-scheme $X$, and consider the corresponding stack $[X/G]$.  Assume that the action is tame or equivalently that the stack $[X/G]$ is tame.  For simplicity, assume that $\bF_q$ contains the $r^{th}$ roots of unity for each $r$ dividing the order of $G$.  We then have $\ii_{\bmu}([X/G]) \cong \ii([X/G])$ the usual inertia stack, and thus
\begin{equation}
\ii_{\bmu}([X/G]) = \bigsqcup_{(g)} [X^g / C(g)].
\end{equation}
Here the union is over conjugacy classes of elements of $G(\bF_{q})$, $C(g)$ denotes the centralizer group scheme of $g$, and $X^g$ denotes the fixed subscheme of $g$.
\end{ex}

Let $\sX$ and $\sY$ be $\bF_q$-stacks.  Let $f$ and $g$ be 1-morphisms from $\sX$ to $\sY$, and let $\phi: f \rightrightarrows g$ be a 2-morphism.  Then composition induces 1-morphisms $\ii_{\bmu}(f)$ and $\ii_{\bmu}(g)$, and a 2-morphism $\ii_{\bmu}(\phi)$ making $\ii_{\bmu}(-)$ into a 2-functor.

Definition \ref{inertia_defn} and the above remarks clearly also apply to $\ov{\bF}_q$-stacks.  The map $\ov{\bF}_q \xrightarrow{\phi} \ov{\bF}_q$ then induces two morphisms on $\ii_{\bmu}(\ov{\sX})$ which agree by the following lemma.

\begin{lem} \label{inertia_compat}
  \hfill
  \begin{enumerate}
    \item There is an equivalence
        $$\ii_{\bmu}(\sX \times_{\bF_q} \ov{\bF}_q) \xrightarrow{\cong} \ii_{\bmu}(\sX) \times_{\bF_q} \ov{\bF}_q.$$
    \item We denote the latter simply by $\ov{\ii_{\bmu}(\sX)}$.  Under the identification above, the following functors are 2-isomorphic:
        $$\ii_{\bmu}(F_{\ov{\sX}})\ \overset{\cong}{\rightrightarrows} \  F_{\ov{\ii_{\bmu}(\sX)}}.$$
  \end{enumerate}
\end{lem}
\begin{proof}
The proof proceeds exactly as in the proof of Lemma \ref{K_compat} appearing in the appendix.
\end{proof}
\begin{rem}
It is important to note that there are several variations of the definition of inertia stack appearing in the literature.  In particular, to work in the most general setting one should consider the \emph{rigidified} inertia stack \cite{AGV1, AGV2}.  However, since our goal is to study the Chen-Ruan cohomology (and hence only degree zero stable maps), Definition \ref{inertia_defn} above is preferred.
\end{rem}
\section{$\ell$-adic Cohomology} \label{l-adic_cohomology}
Let $l$ be coprime to $q$ and denote by $\bZ_l$ (resp. $\bQ_l$) the $\ell$-adic integers (resp. numbers).  Let $\sX$ be a proper, smooth, tame, \DM\ $\bF_q$-stack with projective coarse moduli scheme, and let $\ov{\sX} \coloneq \sX \times_{\bF_q} \ov{\bF}_q$.

\begin{defn} \label{coho_defn}
The \emph{orbifold $\ell$-adic cohomology} of $\ov{\sX}$ is a graded $\bQ_l$-algebra.  As a vector space, it is given by

\begin{alignat}{1}
H^*_{CR}(\ov{\sX}, \bQ_l)\  &\coloneq \ H^*(\ii_{\bmu}(\ov{\sX}), \bQ_l) \notag \\
    &= \ \bigoplus_i \ (\ \underset{n}{\varprojlim}\  H^i(\ii_{\bmu}(\ov{\sX}),\ \bZ / l^n\bZ)\ ) \underset{\bZ_l}{\otimes} \bQ_l. \notag
\end{alignat}
\end{defn}

Any morphism $g: \ov{\sX} \To \ov{\sY}$ induces $\ii_{\bmu}(g): \ii_{\bmu}(\ov{\sX}) \To \ii_{\bmu}(\ov{\sY})$.  Since $\bZ / l^r\bZ$ is a constant sheaf we have $\ii_{\bmu}(g)^*(\bZ / l^r\bZ) \cong \bZ / l^r\bZ$.  Passing to the limit, we obtain a linear map on $H^*(\ii_{\bmu}(\ov{\sX}), \bQ_l)$.  Corresponding push-forward maps on homology groups are then induced by duality.

We now introduce the grading on $H^*(\ii_{\bmu}(\ov{\sX}), \bQ_l)$ using the notion of \emph{age}.  We include the definition at this stage because the definition itself is fairly direct.  The motivation however arises later from the appearance of age in the Riemann-Roch theorem on curves \cite[Thm. 7.2.1]{AGV2} (and hence also appears in the computation of the degree of the virtual fundamental class on the moduli stack of stable maps into $\sX$).  Since one goal of this article is to determine the effect of the Frobenius morphism on this virtual fundamental class, we proceed without assuming our base field is algebraically closed.

For every $r$ fix an embedding $\bmu_r \hookrightarrow \mathbb{G}_m$.  Then a group scheme morphism $\rho: \bmu_r \To \mathbb{G}_m$ over ${\bF}_q$ is determined by an integer $0 \le k \le r-1$ with $\rho(g) = g^k$.  Define
$$\text{age}(\rho) \coloneq \frac{k}{r} \in \bQ.$$
When $(q, r) = 1$, this function extends by linearity to a function on the representation ring of the group scheme.  For any object $((B\bmu_r)_S \xrightarrow{f} \sX)$ in $\ii_{\bmu_r}(\sX)(S)$, each fiber of $f^*T_{\sX}$ over $S$ gives a representation of $\bmu_r$, and we obtain a locally constant function on $\ii_{\bmu_r}(\sX)$.  We obtain a well-defined locally constant function also denoted by  $\text{age}$:
$$\text{age}: \ii_{\bmu}(\sX) \To \bQ.$$
Finally, we may pull this function back to a function on $\ii_{\bmu}(\ov{\sX})$, and this induces a function $\text{age}: H^*_{CR}(\ov{\sX}, \bQ_l) \To \bQ$.  The various uses of this notation will be clear from the context.  The grading is then given by
\begin{equation} \label{grading}
H^*_{CR}(\ov{\sX}, \bQ_l) = \bigoplus_{i \in \bQ} H_{CR}^{i}(\ov{\sX}, \bQ_l)
\end{equation}
where
\begin{equation} \label{grading2}
H_{CR}^{i}(\ov{\sX}, \bQ_l) \coloneq \bigoplus_{a + 2b = i} H^a(\text{age}^{-1}(b), \bQ_l).
\end{equation}

Note that the choices of embeddings $\bmu_r \hookrightarrow \mathbb{G}_m$ will change the age function on $\ii_{\bmu}(\ov{\sX})$.
\begin{rem} \label{inertia_remark}
Note that in general an automorphism of $\ii_{\bmu}(\sX)$ \emph{over} $\sX$ will not preserve the age function.  An essential example is the involution $i: \ii_{\bmu}(\sX) \To \ii_{\bmu}(\sX)$ given by precomposing each $(B\bmu_r)_S \To \sX$ with the automorphism $(B\bmu_r)_S \To (B\bmu_r)_S$ induced by $g \mapsto g^{-1}: \bmu_r \To \bmu_r$.  On the other hand, any automorphism of $\ii_{\bmu}(\sX)$ \emph{induced by an automorphism of $\sX$} does preserve the age.  An essential example is the arithmetic Frobenius on $\ov{\sX}$.
\end{rem}
\begin{rem} \label{Yasuda}
In \cite{Yas}, Yasuda defines the additive Chen-Ruan cohomology as in Definition \ref{coho_defn} above with the following exceptions.  If we define a function on $\ii_{\bmu}(\ov{\sX})$ by $sht := age \circ i$ where $i$ is as in Remark \ref{inertia_remark}, then Yasuda replaces the right side of (\ref{grading2}) with
$$\bigoplus_{a + 2b = i} H^a(\text{sht}^{-1}(b), \bQ_l(-b)).$$  (Strictly speaking, Yasuda also works with the cohomology of coarse moduli space.)  The Tate twist $\bQ_l(-b)$ changes the weight as a Galois representation, and one motivation for this arises from Yasuda's proof (\cite[Cor. 4.9]{Yas}) of an additive version of the Crepant Resolution Conjecture \cite{BG, Ruan1, Ruan2}.  In the current paper, another motivation is found via the Galois action on the \emph{ring} structure of Chen-Ruan cohomology defined in the section.
\end{rem}
\section{Ring Structure}
In this section we gather the results in \cite{BF}, \cite{AV1}, and \cite{AGV2} needed to define the ring structure.  One subtlety is the cycle map from Chow groups to $\ell$-adic cohomology which requires the moduli stack to be smooth (see Remark \ref{cycle}).

\subsection*{Introduction}
The ring structure we shall impose on $H^*_{CR}(\ov{\sX}, \bQ_l)$ is motivated by quantum cohomology, and so we give a brief description.  For simplicity consider a smooth projective scheme $Y$ over $\bC$.  The quantum cohomology of $Y$ is a deformation of $H^*(Y, \bC)$.  If $H_2^+(Y,\bZ) \subset H_2(Y, \bZ)$ denotes the classes generated by effective curves in $Y$, then the parameter space of the deformation is given by the semigroup algebra $\bQ[H_2^+(Y, \bZ)]$.  The product in this deformed ring requires integrals over the moduli stack of curve in $Y$ (more precisely over a compactification by stable maps due to Kontsevich \cite{Ko}).  Finally, the original ring $H^*(Y, \bC)$ is recovered by setting the deformation parameters to zero (i.e. by only considering the moduli stack of constant stable maps).  However, when $Y$ is replaced with a stack $\sY$, this limit of quantum cohomology does not agree with $H^*(\sY, \bC)$.  The new ring obtained is called the \emph{Chen-Ruan cohomology ring}.  In what follows we define the ring structure directly using the complex case as motivation.

\subsection*{Stable maps}Let $\sX$ be a proper, smooth, tame, \DM\ stack over $\bF_q$ with projective coarse moduli scheme $X$, and let $\ov{\sX} \coloneq \sX \times_{\bF_q} \ov{\bF}_q$.  We define the moduli stack of stable maps into $\sX$ as constructed in \cite{AV1}.

\begin{defn}
A \emph{balanced twisted stable n-pointed map} $(\sC \xrightarrow{f} \sX_S, \{\Sigma_i\})$ is a commutative diagram of $\bF_q$-stacks
\begin{diagram}
\sqcup_{i=1}^n \Sigma_i & \rInto & \sC  & \rTo^{f}            & \sX_S          \\
                        &        & \dTo^{\pi} &                     &  \dTo        \\
                        &        & C    & \rTo^{|f|}       & X_S
\end{diagram}
where
\begin{enumerate}
  \item $\sC$ is a proper \DM\ stack with coarse moduli space $C$
  \item $(C, \{\pi(\Sigma_i)\})$ is an n-pointed nodal curve
  \item Over the node $\{xy = 0\}$ of $C$, $\sC$ has \'{e}tale chart $$[\{xy = 0\} / (\bmu_r)_{\bF_q}]$$ where the action is given by $(x,y) \mapsto (\xi u, \xi^{-1} v)$
  \item Over a marked point $\pi(\Sigma_i)$ of $C$, $\sC$ has \'{e}tale chart $$[\mathbb{A}^1 / (\bmu_r)_{\bF_q}]$$ where the action is given by $u \mapsto \xi u$ and $\Sigma_i$ is the substack defined by $u = 0$
  \item $\pi$ is an isomorphism away from the markings and nodes
  \item $f$ is representable with $|f|$ the induced map on coarse moduli spaces
  \item $|f|$ is stable in the sense of Kontsevich \cite{Ko}.
\end{enumerate}
\end{defn}

We shall refer to the above merely as \emph{stable maps}.

For any $d,g \ge 0$, we say that $(\sC \xrightarrow{f} \sX, \{\Sigma_i\})$ has degree $d$ and genus $g$ if $|f|$ does.  After appropriately defining such objects over an arbitrary $\bF_q$-scheme, one obtains the stack $\sK_{g,n}(\sX, d)$, a proper stack of finite type over $\bF_q$ with projective coarse moduli space \cite[Theorem 1.4.1]{AV1}.

Our purposes only require the case when $g = 0$, $d = 0$, and $n = 3$.  We shall denote $\sK_{0,3}(\sX, 0)$ simply by $\sK(\sX)$, and we have the following additional properties:

\begin{lem} \label{K_smooth}
$\sK(\sX) = \sK_{0,3}(\sX, 0)$ is a smooth \DM\ stack over $\bF_q$.
\end{lem}
\begin{proof}
See the appendix.
\end{proof}
\begin{lem-def} \label{eval}
There exist evaluation morphisms over $\bF_q$ denoted $e_i$ for $i = 1,2,3:$
$$\sK(\sX) \xrightarrow{e_i} \ii_{\bmu}(\sX)$$
which applies to objects over $\bF_q$ as
$$(\sC \xrightarrow{f} \sX, \{\Sigma_i\}) \mapsto (\Sigma_i \xrightarrow{f|_{\Sigma_i}} \sX).$$
\end{lem-def}
\begin{proof}
The proof of \cite[Lemma 6.2.1]{AGV1} carries over to this context without change.
\end{proof}
\begin{rem}
The subtlety of this lemma lies in the definition of stable maps over more general base schemes.
\end{rem}

\subsection*{Virtual Classes}Finally, we consider integrating cohomology classes on this moduli space, and for this purpose we a need a fundamental class.  However, even though $\sK(\sX)$ is smooth, the natural vector spaces holding obstructions to deforming stable maps may still be non-trivial.  In \cite{BF}, the authors define an \emph{obstruction theory} to describe this phenomenon.  They proceed to construct a \emph{virtual fundamental class} in the Chow group of the moduli stack as a replacement of the usual fundamental class.  In the case of interest in this article, these constructions have the following concrete description.

Let
\begin{equation} \label{universal_diagram}
\begin{diagram}
\mathcal{U} & \rTo^{f} & \sX \\
\dTo_{\pi} & & \\
\sK(\sX) & & \\
\end{diagram}
\end{equation}
be the universal curve and universal stable map to $\sX$.  Then we have the following lemma.

\begin{lem} \label{obs_theory}
    \hfill
    \begin{enumerate}
        \item The natural map $$(R^{\bullet}\pi_* f^*T_{\sX})^{\vee} \ \xrightarrow{\phi} \ \Omega^1_{\sK(\sX)/ \sM_{0,3}^{tw}}$$ is a perfect relative obstruction theory with virtual dimension (denoted vdim) given by the locally constant function $$vdim  = dim \sX - \text{age}\circ e_1 - \text{age}\circ e_2 - \text{age}\circ e_3.$$
        \item $R^1\pi_* f^*T_{\sX}$ is locally free (denote the locally constant rank by $r$), and the virtual fundamental class (denoted $[\sK(\sX)]^{vir}$) in $A_{vdim}(\sK(\sX))_{\bQ}$ induced by $\phi$ is $$[\sK(\sX)]^{vir} = c_{r}(R^1\pi_*f^*T_{\sX}).$$
    \end{enumerate}
\end{lem}
\begin{proof}
See the appendix.
\end{proof}
\begin{rem}
The previous lemma will also hold if $\sX$ is replaced by $\ov{\sX}$.  However, the statement over the non-algebraically closed field is essential in the proof of Proposition \ref{main} (specifically  Lemma \ref{F_on_vir_class}) where we compute the action of the Frobenius on $[\sK(\ov{\sX})]^{vir}$ in the Chow group.
\end{rem}
\begin{rem} \label{cycle}
Since $\sK(\ov{\sX})$ is smooth, one can construct the cycle map $$A_*(\sK(\ov{\sX}))_{\bQ} \otimes_{\bQ} \bQ_l \xrightarrow{cl} H_*(\sK(\ov{\sX}), \bQ_l) = H^*(\sK(\ov{\sX}), \bQ_l)$$ by proceeding as in \cite[VI. Section 9]{Milne} using the Gysin sequence in \cite[Cor. 2.1.3]{Be-Lef} and a slight refinement of the long exact sequence in Section 2.1 of [ibid.].  We denote the image of $[\sK(\ov{\sX})]^{vir}$ under $cl$ also by $[\sK(\ov{\sX})]^{vir}$.
\end{rem}

\subsection*{Ring structure}
The ring structure can now be constructed formally on $H^*_{CR}(\ov{\sX}, \bQ_l)$ just as in the original formulation \cite{CR2}.  For $\alpha, \beta \in H^*_{CR}(\ov{\sX}, \bQ_l)$ define $$\alpha \star \beta \coloneq i_*(e_3)_*(e_1^*\alpha \cup e_2^* \beta \cap [\sK(\ov{\sX})]^{vir})$$ where $i: \ii_{\bmu}(\ov{\sX}) \To \ii_{\bmu}(\ov{\sX})$ is the morphism induced by the isomorphisms $\lambda \mapsto \lambda^{-1}: \bmu_r \To \bmu_r$.
\begin{prop} \label{assoc}
The operation $\star$ makes $H^*_{CR}(\ov{\sX}, \bQ_l)$ into a (graded) commutative, associative ring with unity.
\end{prop}
\begin{proof}
The description of the boundary strata of $\sK_{0,4}(\ov{\sX}, 0)$, and the proofs in Sections 5 and 6 of \cite{AGV2} apply to this context as well.  Note that since we are only concerned with degree zero maps, there is no need to decompose $\sK_{0,n}(\ov{\sX})$ via curve classes on the coarse moduli space of $\ov{\sX}$.
\end{proof}

\begin{ex} \label{example_ring}
Suppose $b$ is coprime to $q$ and $\bG_m = \text{Spec }(\bF_q[t, \frac{1}{t}])$ acts on $\bA^2$ with weights $1$ and $b$.  Assume for simplicity that $b$ is prime and $F_q$ contains the $b^{th}$ roots of unity.  The stack $\sX \coloneq [(\bA^2 \backslash \{0\})/ \bG_m]$ has \'{e}tale neighborhoods
\begin{alignat}{1}
[\bA^1/\bmu_b] &\To \sX \notag \\
\bA^1 &\To \sX \notag
\end{alignat}
whose induced map on coarse moduli spaces form a cover of $\bP^1$, the coarse moduli scheme of $\sX$.  Here the $\bmu_b$ action is given by the character $\lambda \mapsto \lambda^b: \bmu_b \To \bG_m$.  We may then use Example \ref{inertia_example} above to compute
$$\ii_{\bmu}(\sX) = \sX\ \bigsqcup\ \sqcup_{i=1}^{b-1} B\bmu_b.$$
Let $A$ denote the copy of $\ov{B\bmu_b} \subset \ii_{\bmu}(\ov{\sX})$ with age equal to $\frac{1}{b}$.  Denote also by $A$ a generator of $H^0(A, \bQ_l) \subset H^0(\ii_{\bmu}(\ov{\sX}), \bQ_l)$.  Then one has
$$H^*_{CR}(\ov{\sX}, \bQ_l) \cong \bQ_l[A]/\langle A^{b+1} \rangle.$$
In particular, $A^b$ generates $H^2(\ov{\sX}, \bQ_l) \subset H^2(\ii_{\bmu}(\ov{\sX}), \bQ_l)$ and has age equal to $0$.
\end{ex}
\section{Frobenius Actions} \label{Frobenius_Actions}
Let $\sX$, $\ov{\sX}$ be as in the last section.  The arithmetic Frobenius morphism $F_{\ov{\sX}}$ on $\ov{\sX}$ induces the morphism $\ii_{\bmu}(F_{\ov{\sX}}) = F_{\ii_{\bmu}(\ov{\sX})}$ on $\ii_{\bmu}(\ov{\sX})$.  Thus we have an induced map on Chen-Ruan cohomology groups which we denote by
$$\ii_{\bmu}(F_{\ov{\sX}})^* : H^*_{CR}(\ov{\sX}, \bQ_l) \To H^*_{CR}(\ov{\sX}, \bQ_l).$$

Consider Example \ref{example_ring} above.  $\ii_{\bmu}(F_{\ov{\sX}})^*$ acts on the orbifold $\ell$-adic cohomology by fixing $A^1, \cdots, A^{b-1}$ and sending $A^b$ to $q^{-1}A^b$.  Thus $\ii_{\bmu}(F_{\ov{\sX}})^*$ does not preserve the ring structure.  However by composing $\ii_{\bmu}(F_{\ov{\sX}})^*$ with the morphism sending $A^i \mapsto q^{-\text{age}(A^i)}A^i$ we indeed obtain a homomorphism of rings sending $A^i \mapsto q^{-\frac{i}{b}} A^i$. Proposition \ref{main} from the introduction shows that this phenomenon holds in general.

\begin{no-num}{Proposition \ref{main}}
The \emph{orbifold Frobenius morphism} given by
\begin{alignat}{1}
F_{\ov{\sX}, q, orb}: \ H^*_{CR}(\ov{\sX}, \bQ_l) &\longrightarrow H^*_{CR}(\ov{\sX}, \bQ_l) \notag \\
                                        \alpha       &\longmapsto q^{-\text{age}(\alpha)} \cdot \ii_{\bmu}(F_{\ov{\sX}, q})^*(\alpha) \notag
\end{alignat}
is a homomorphism of graded rings.
\end{no-num}
When no confusion arises, the subscripts $\ov{\sX}$ and $q$ on $F$ may be dropped.
\begin{rem} \label{Yasuda2}
Note that instead of twisting the natural map $\ii_{\bmu}(F_{\ov{\sX}})^*$ to obtain $F_{orb}$, one could twist the coefficients of the Chen-Ruan cohomology.  For instance, if we redefined the Chen-Ruan cohomology groups by replacing the right side of (\ref{grading2}) with
$$\bigoplus_{a + 2b = i} H^a(\text{age}^{-1}(b), \bQ_l(-b)),$$
then the natural map action of $\ii_{\bmu}(F_{\ov{\sX}})$ on these new groups would be a ring isomorphism.
\end{rem}

Before proving the proposition we make a few observations.  First we consider the map on $\sK(\ov{\sX})$ induced by $F_{\ov{\sX}}$.  For any $\ov{\bF}_q$-scheme $S$, and any object
\begin{diagram}
\sqcup_{i=1}^3 \Sigma_i & \rInto & \sC  & \rTo^{f_S}            & \ov{\sX}          \\
                        &        & \dTo &                     &  \dTo        \\
                        &        & S    & \rTo^{\pi}      & \text{Spec}(\ov{\bF}_q)
\end{diagram}
of $\sK(\ov{\sX})(S)$, the commutative diagram
\begin{diagram}
\sqcup_{i=1}^3 \Sigma_i & \rInto & \sC  & \rTo^{F_{\ov{\sX}} \circ f_S}            & \ov{\sX}          \\
                        &        & \dTo &                     &  \dTo        \\
                        &        & S    & \rTo^{F_{\ov{\bF}_q} \circ \pi}      & \text{Spec}(\ov{\bF}_q)
\end{diagram}
is clearly an object of $\sK(\ov{\sX})$ once we check the stability condition.  However, this is easy since $f_S$ has degree zero.  The obvious map on morphisms then determines a functor we denote by $\sK(F_{\ov{\sX}})$.  Note that $\sK(F_{\ov{\sX}})$ covers the Frobenius map on $\text{Spec}(\ov{\bF}_q)$.

The following lemma compares $\sK(F_{\ov{\sX}})$ with the Frobenius morphism $F_{\sK(\ov{\sX})}$.
\begin{lem} \label{K_compat}
  \hfill
  \begin{enumerate}
    \item There is an equivalence
        $$\sK(\sX \times_{\bF_q} \ov{\bF}_q) \xrightarrow{\cong} \sK(\sX) \times_{\bF_q} \ov{\bF}_q.$$
    \item We denote the latter simply by $\ov{\sK(\sX)}$.  Under the identification above, the following functors are 2-isomorphic:
        $$\sK(F_{\ov{\sX}})\ \overset{\cong}{\rightrightarrows} \  F_{\ov{\sK(\sX)}}.$$
  \end{enumerate}
\end{lem}
\begin{proof}
See the appendix.
\end{proof}

These identifications allow us to determine the image of the virtual fundamental class under $\sK(F_{\ov{\sX}})$ in homology.
\begin{lem} \label{F_on_vir_class}
$\sK(F_{\ov{\sX}})_* [\sK(\ov{\sX})]^{vir} = q^{-vdim}[\sK(\ov{\sX})]^{vir}$
\end{lem}
\begin{proof}
See the appendix.
\end{proof}
\begin{proof}[Proof of Proposition \ref{main}]
For simplicity, denote $F_{orb} \coloneq F_{\ov{\sX}, q, orb}$, $\ii(-) \coloneq \ii_{\bmu}(-)$, and $F \coloneq F_{\ov{\sX}, q}$.  Let $\ii(\ov{\sX}) = \bigsqcup_j \ii(\ov{\sX})_j$ be a decomposition into connected components.  Define $a_j \coloneq \text{age}(\ii(\ov{\sX})_j)$ and $\hat{a}_j \coloneq \text{age}(i(\ii(\ov{\sX})_j))$ where $i: \ii_{\bmu}(\sX) \To \ii_{\bmu}(\sX)$ is the isomorphism appearing in the definition of $\star$ above.  To prove the proposition it suffices to check that
\begin{equation} \label{compat}
    F_{orb}(\alpha_1) \star F_{orb}(\alpha_2) = F_{orb}(\alpha_1 \star \alpha_2)
\end{equation}
for any $\alpha_k \in H^*(\ii(\ov{\sX})_{j(\alpha_k)}, \bQ_l) \subset H^*(\ii(\ov{\sX}), \bQ_l)$ ($k = 1,2$).  Apply $\ii(F)^* \ii(F)_*$ to the left side of (\ref{compat}), where Poincar\'{e} duality isomorphisms have been suppressed.  By the projection formula and since $\sK(F)$ and $\ii(F)$ commute, we then have
\begin{alignat}{1} \label{eqn1}
\ii(F)^* \ii(F)_* (&F_{orb}(\alpha_1) \star F_{orb}(\alpha_2)) = \notag \\
    &= \ii(F)^* i_* (e_3)_* \sK(F)_* (q^{-a_1 - a_2} \sK(F)^* e_1^* \alpha_1 \cup \sK(F)^* e_2^* \alpha_2 \cap [\sK(\ov{\sX})]^{vir}) \notag \\
    &= \ii(F)^* i_* (e_3)_*(q^{-a_1 - a_2} e_1^* \alpha_1 \cup e_2^* \alpha_2 \cap \sK(F)_* [\sK(\ov{\sX})]^{vir}) \notag \\
    &= \ii(F)^* i_* (e_3)_*(q^{-a_1 - a_2 - vdim} e_1^* \alpha_1 \cup e_2^* \alpha_2 \cap [\sK(\ov{\sX})]^{vir}).
\end{alignat}
Now note that the operator $\ii(F)^* \ii(F)_*$ decomposes according to the following lemma whose proof is left to the reader:
\begin{lem}
Let $1_j$ denote the denote both the fundamental class of $\ii(\ov{\sX})_j$ and the operator given by taking cup product with $1_j$.  Then $$\ii(F)^* \ii(F)_* = \sum_j q^{-\text{dim} \ii(\ov{\sX})_j} 1_j.$$
\end{lem}
Thus (\ref{eqn1}) implies
\begin{alignat}{1} \label{eqn2}
  F_{orb}(\alpha_1) &\star F_{orb}(\alpha_2) = \notag \\
  &= \sum_j 1_j\ q^{\text{dim} \ii(\ov{\sX})_j} \ii(F)^*i_* (e_3)_*(q^{-a_1 - a_2 - \text{vdim}} e_1^* \alpha_1 \cup e_2^* \alpha_2 \cap [\sK(\ov{\sX})]^{vir}) \notag \\
  &= \sum_j 1_j\ F_{orb} i_* (e_3)_*(q^{-a_1 - a_2 - \text{vdim} + \text{dim} \ii(\ov{\sX})_j + a_3} e_1^* \alpha_1 \cup e_2^* \alpha_2 \cap [\sK(\ov{\sX})]^{vir}).
\end{alignat}
Now for each $j$, the $1_j$ in (\ref{eqn2}) restricts the class to $\ii(\ov{\sX})_j$.  This contribution will not be changed if we replace $[\sK(\ov{\sX})]^{vir}$ by its restriction to $e_1^{-1}\ii(\ov{\sX})_{j_1} \bigcap \ e_2^{-1}\ii(\ov{\sX})_{j_2} \bigcap\ (i \circ e_3)^{-1}\ii(\ov{\sX}_j)$ where $j_1$ (resp. $j_2$) is the index of the component of $\ii(\ov{\sX})$ supporting $\alpha_1$ (resp. $\alpha_2$).  On this locus of $\sK(\ov{\sX})$, vdim is constant and equal to $\text{dim}\ \ov{\sX} - a_1 - a_2 - \hat{a}_3$.

Since for each $k$, $$a_k + \hat{a}_k = \text{dim}\ \ov{\sX} - \text{dim} \ii(\ov{\sX})_{j(a_k)},$$ the exponent of $q$ in (\ref{eqn2}) is zero.  Thus (\ref{eqn2}) implies
\begin{alignat}{1}
F_{orb}(\alpha_1) \star F_{orb}(\alpha_2)
    &= \sum_j 1_j\ F_{orb} i_* (e_3)_*(e_1^* \alpha_1 \cup e_2^* \alpha_2 \cap [\sK(\ov{\sX})]^{vir}) \notag \\
    &= F_{orb}(\alpha_1 \star \alpha_2), \notag
\end{alignat}
and the proposition is proved.
\end{proof}
\section{Orbifold Zeta Functions} \label{zeta}
Let $\sX$, $\ov{\sX}$ be as in the last section.  For simplicity assume in addition that $\sX$ is \emph{Gorenstein} condition:  $\sX$ has generically trivial isotropy, and for any geometric point $\xi \in \sX(\ov{\bF}_q)$, the representation $$\text{Aut}(\xi) \To \text{GL}(T_{\ov{\sX}}|_{\xi})$$ has determinant 1.  The latter condition is equivalent to $\text{age}: \ii_{\bmu}(\sX) \To \bQ$ taking integer values.  Thus $H^*_{CR}(\ov{\sX}, \bQ_l)$ is $\bZ$-graded.

For a linear map $F: V \To V$ on a $\bZ$-graded vector space we write $V = \oplus V_i$, $F = \oplus F_i$ and we denote
\begin{alignat}{1}
\text{det}\ (F\ |\ V)\ &=\ \prod_i \text{det}\ (F_i\ |\ V_i)^{(-1)^{i+1}} \notag \\
\text{Tr}\ (F\ |\ V)\ &=\ \sum_i (-1)^i\ \text{Tr}\ (F_i\ |\ V_i). \notag
\end{alignat}

\begin{defn} \label{zeta_defn}
The \emph{orbifold cohomological zeta function} is given by
\begin{alignat}{1}
Z_{H^*_{CR}}(\sX, t)\ &\coloneq\ \text{det}\ (1 - F_{orb}\ t\ |\ H^*_{CR}(\ov{\sX}, \bQ_l)) \notag \\
                      &=\ \text{exp}\ (\sum_{r=1}^{\infty} \text{Tr}(F_{orb}^{r}\ |\ H^*_{CR}(\ov{\sX}, \bQ_l)) \frac{t^r}{r}). \notag
\end{alignat}
\end{defn}

One obtains a trace formula for $Z_{H^*_{CR}}(\sX, t)$ by applying the Lefschetz Trace Theorem of \cite{Be-Lef} to $\ii_{\bmu}(\ov{\sX})$.  For a $\bF_q$-scheme $S$, let $[\ii_{\bmu}(\sX)(S)]$ denote the set of isomorphism classes of the groupoid $\ii_{\bmu}(\sX)(S)$.  For $\xi \in [\ii_{\bmu}(\sX)(S)]$, let $\text{Aut}(\xi)$ denote the automorphism group of any representative of $\xi$.  Finally, let $\ii_{\bmu}(\sX) = \bigsqcup_i \ii_{\bmu}(\sX)_i$ be a decomposition into connected components so that $\text{age}$ and dimension ($\text{dim}$) are constant on each $\ii_{\bmu}(\sX)_i$.  Then \cite[Theorem 3.1.2]{Be-Lef} yields
\begin{alignat}{1} \label{trace_formula}
\text{Tr}\ (F_{orb}\ |\ H^*_{CR}(\ov{\sX}, \bQ_l))\
    &=\ \sum_i \text{Tr}\ (F_{orb}|_{\ii_{\bmu}(\ov{\sX})_i}\ |\ H^*(\ii_{\bmu}(\ov{\sX})_i, \bQ_l)) \notag \\
    &=\ \sum_i q^{-\text{age}(\ii_{\bmu}(\ov{\sX})_i)} \text{Tr}\ (F^*_{\ii_{\bmu}(\ov{\sX})_i}\ |\ H^*(\ii_{\bmu}(\ov{\sX})_i, \bQ_l)) \notag \\
    &=\ \sum_i q^{-\text{age}(\ii_{\bmu}(\ov{\sX})_i)} \sum_{\xi \in [\ii_{\bmu}(\sX)_i(\bF_q)]} \frac{q^{-\text{dim}(\xi)}}{\# \text{Aut}(\xi)} \notag \\
    &=\ \sum_{\xi \in [\ii_{\bmu}(\sX)(\bF_q)]} \frac{q^{-\text{age}(\xi)-\text{dim}(\xi)}}{\# \text{Aut}(\xi)}.
\end{alignat}

\begin{rem}
One is led to claim that the trace of $F_{orb}$ counts objects of \\$\bigcup_r \text{\underline{Hom}}^{\text{rep}}_{\ \bF_q}(B\bmu_r, \sX)(\bF_q)$ counted with weights by the age and dimension.  It is natural to ask if the trace of $F_{orb}$ counts some natural objects on $\sX$ without weights.
\end{rem}

Recall that $F_{orb} = F_{q, orb}$ depends on the base field.  Since $F^r_{q, orb} = F_{q^r, orb}$, equation (\ref{trace_formula}) yields a formula for the trace of each iterate of $F_{q, orb}$.  We then obtain the following analog of (\ref{intro_eqn_1}):
\begin{equation}
\text{det}\ (1 - F_{orb}\ t\ |\ H^*_{CR}(\ov{\sX}, \bQ_l))\ =\ \text{exp}\ \sum_{r = 1}^{\infty}\ ( \sum_{\xi \in [\ii_{\bmu}(\sX)(\bF_{q^r})]} \frac{(q^r)^{-\text{age}(\xi)-\text{dim}(\xi)}}{\# \text{Aut}(\xi)} ) \frac{t^r}{r}.
\end{equation}

We attempt to further interpret the arithmetic information contained by $Z_{H^*_{CR}}(\sX, t)$ motivated by the Crepant Resolution Conjecture \cite{BG, Ruan1, Ruan2}.  A theorem of Yasuda \cite[Cor. 4.9]{Yas} gives the following result.  Two smooth, proper stacks $\sX_1$ and $\sX_2$ are \emph{K-equivalent} if there exists smooth, proper stack $\sY$ and proper, tame and birational maps $\pi_1$ and $\pi_2$
\begin{diagram}
        &               & \sY   &               & \\
        & \ldTo^{\pi_1} &       & \rdTo^{\pi_2} & \\
\sX_1   &               &       &               & \sX_2 \\
\end{diagram}
such that $\pi_1^*K_{\sX_1} \cong \pi_2^*K_{\sX_2}$ where $K_{\sX_i}$ is the canonical line bundle on $\sX_i$.
\begin{thm-no-num}
If $\sX_1$ and $\sX_2$ are K-equivalent, proper, smooth, tame \DM\ stacks (recall also the Gorenstein assumption made throughout this section), then
$$Z_{H^*_{CR}}(\sX_1, t) = Z_{H^*_{CR}}(\sX_2, t).$$
\end{thm-no-num}
\begin{proof}
We simply show that the orbifold zeta function defined above agrees with the natural zeta function built from the Galois representation in \cite[Defn. 4.6]{Yas} for which the statement holds \cite[Cor. 4.9]{Yas}.  Following Remarks \ref{Yasuda} and \ref{Yasuda2} above, it suffices to show
\begin{alignat}{1}
\text{Tr}\ (\ \ii(F)^* \ |\ \bigoplus_{a + 2b = j} &H^a(\text{age}^{-1}(b), \bQ_l(-b))\ ) = \notag \\
    &\text{Tr}\ (\ \ii(F)^* \ |\ \bigoplus_{a + 2b = j} H^a((\text{age} \circ i)^{-1}(b), \bQ_l(-b))\ ) \notag
\end{alignat}
where $i: \ii(\sX) \To \ii(\sX)$ is the involution.  But this follows from $i \circ \ii(F) = \ii(F) \circ i$ which one can easily show.
\end{proof}
In particular, we see that the orbifold zeta function carries the arithmetic information of any crepant resolution (when one exists) of the coarse moduli scheme.
\begin{cor-no-num} \label{Yasuda_corollary}
Let $\sX$ be a proper, smooth, tame \DM\ stack satisfying the hard Lefschetz condition with trivial generic stabilizer.  Suppose $Y \To X$ is a crepant resolution of the coarse moduli scheme $X$ of $\sX$, then
$$Z_{H^*_{CR}}(\sX, t) = Z(Y,t)$$
where $Z(Y,t)$ is the classical zeta function.
\end{cor-no-num}
\begin{rem}
Yasuda also proves an analog of \cite[Cor. 4.9]{Yas} over the complex numbers, and this result was obtained independently by Lupercio and Poddar \cite{LP}.
\end{rem}
\begin{rem}
It is natural to associate to an orbifold, the zeta function of any crepant resolution (when one exists) of the coarse moduli space (see for example \cite[page 9]{W}).  The above corollary then shows that this definition agrees with Definition \ref{zeta_defn} above in the special case when such a crepant resolution exists.
\end{rem}

\begin{ex}
Suppose 2 is coprime to $q$, and suppose $\bG_m = \text{Spec }(\bF_q[t, \frac{1}{t}])$ acts on $\bA^3$ with weights $1$, $1$, and $2$.  The stack $\sX \coloneq [\bA^3 \backslash \{0\}/ \bG_m]$ has \'{e}tale neighborhood $[\bA^2/ \bmu_2] \To \sX$ where the action of $\bmu_2 \cong \bZ/2\bZ$ is given by the direct sum of two copies of the non-trivial character of $\bmu_2$.  If $$\sX \xrightarrow{\pi_1} |\sX|$$ denotes the morphism to the coarse moduli scheme, then $|\sX|$ is the projective closure of $\bA^2/ \bmu_2 \cong \text{Spec }(\bF_q[x,y,z]/\langle xy-z^2 \rangle)$.  The canonical sheaf $K_{|\sX|}$ is locally free and $\pi_1^* K_{|\sX|} \cong K_{\sX}$.  Furthermore, the blow-down map $$Y \coloneq \bP(\sO_{\bP^1} \oplus \sO_{\bP^1}(1)) \xrightarrow{\pi_2} |\sX|$$ is a resolution of singularities with $\pi_2^* K_{|\sX|} \cong K_Y$.  Thus the fibered product $Y \times_{|\sX|} \sX$ induces a $K$-equivalence between $Y$ and $\sX$. On each space, the cohomology is generated by algebraic classes and so the Frobenius action is easily computed.

Since $Y$ is a scheme, $\ii_{\bmu}(Y) = Y$, $H^*_{CR}(\ov{Y}, \bQ_l) = H^*(\ov{Y}, \bQ_l)$, and $F_{orb} = F_Y^*$ is the map on cohomology induced by the usual arithmetic Frobenius morphism.  The cohomology of $Y$ is well-known, and we have $$Z_{H^*_{CR}}(Y, t) = \frac{1}{(1-t)(1-q^{-1}t)^2(1-q^{-2}t)}.$$

For the cohomology of $\sX$, we proceed as in Example \ref{example_ring} obtaining
$$\ii_{\bmu}(\sX) = \sX\ \sqcup B\bmu_2.$$
The substack $B\bmu_2$ has age 1 and $H^0(\ov{B\bmu_2}, \bQ_l)$ is fixed by $F_{orb}$.  Thus $B\bmu_2$ contributes a factor of $\frac{1}{1-q^{-1}t}$ to $Z_{H^*_{CR}}(\sX, t)$.  The substack $\sX \subset \ii_{\bmu}(\sX)$ has age 0 and the action of $F_{orb}$ on $H^*(\ov{\sX}, \bQ_l)$ agrees with the action of $F_{|\ov{\sX}|}$ on $H^*(|\ov{\sX}|, \bQ_l)$.  Thus $\sX$ contributes a factor of $\frac{1}{(1-t)(1-q^{-1}t)(1-q^{-2}t)}$ to $Z_{H^*_{CR}}(\sX, t)$, and we have
$$Z_{H^*_{CR}}(\sX, t) = \frac{1}{(1-t)(1-q^{-1}t)^2(1-q^{-2}t)}.$$
\end{ex}
\section{Appendix}
Here we collect the proofs of several lemmas used above.
\begin{no-num}{Lemma \ref{K_smooth}}
$\sK(\sX) = \sK_{0,3}(\sX, 0)$ is a smooth \DM\ stack over $\bF_q$.
\end{no-num}
\begin{proof}
First, the moduli of stable maps to the coarse moduli scheme $X$ is given by $\sK_{0,3}(X, 0) \cong \ov{\sM}_{0,3} \times X \cong \text{Spec}(\bF_q) \times X$, and hence is a \DM\ stack over $\bF_q$.  Thus by \cite[Theorem 1.4.1]{AV1}, $\sK(\sX)$ is a \DM\ stack as well.

To see that $\sK(\sX)$ is smooth, it is sufficient to prove smoothness of $\sK(\ov{\sX})$.  let $\sK(\ov{\sX}) \xrightarrow{p} \sM_{0,3}^{tw}$ be the forgetful functor to the $\ov{\bF}_q$-stack of (not necessarily stable) genus zero twisted curves with three marked points (see \cite{Ol}).  It follows from Remark 1.10 in [ibid.] that $\sM_{0,3}^{tw}$ is smooth over $\ov{\bF}_q$, thus it suffices to show that $p$ is smooth.  By Proposition 17.10 and Corollary 17.9.2 in \cite{LMB}, $p$\ is smooth if and only if $\Omega^1_{\sK(\ov{\sX})/ \sM_{0,3}^{tw}}$ is locally free of finite rank.  For this is suffices to show $R^0 \pi_* f^* T_{\ov{\sX}}$ is locally free of finite rank where $\pi$ and $f$ are the universal curve and stable map respectively as in (\ref{universal_diagram}).

Let $\text{Spec}(\ov{\bF}_q) \xrightarrow{p} \sK(\ov{\sX})$ be a geometric point corresponding to the stable map $(\sC \xrightarrow{g} \ov{\sX}, \{\Sigma_i\})$.  Then to show $R^0 \pi_* f^* T_{\ov{\sX}}$ is locally-free it suffices to show that $\text{dim} H^0(\sC, g^*T_{\ov{\sX}})$ is locally constant as $p$ varies.

For a tuple $\underline{b} = (b_1, b_2, b_3)$ of positive integers, let $\sM_{0,3}^{tw}(\underline{b}) \subset \sM_{0,3}^{tw}$ denote the locus of curves with isotropy group $\bmu_{b_i}$ at the $i^{th}$ marked point.  The decomposition $\sM_{0,3}^{tw} = \bigsqcup_{\underline{b}} \sM_{0,3}^{tw}(\underline{b})$ \cite[Section 5.4]{Ol} induces a decomposition $$\sK(\ov{\sX}) = \bigsqcup_{\underline{b}} p^{-1}(\sM^{tw}_{0,3}(\underline{b})).$$
Note that since $\sK(\ov{\sX})$ consists of degree zero maps, the image of $p$ is contained in the locus of smooth curves with stable coarse moduli space.  Fix $\underline{b}$ for which $p^{-1}(\sM^{tw}_{0,3}(\underline{b}))$ is non-empty and let $\sC$ be the \emph{unique} domain curve of maps in $p^{-1}(\sM^{tw}_{0,3}(\underline{b}))$.  Then we have a morphism of $\ov{\bF}_q$-stacks
$$p^{-1}(\sM^{tw}_{0,3}(\underline{b})) \xrightarrow{\Phi} \text{Pic}(\sC)$$
sending $(\sC \times S \xrightarrow{f} \ov{\sX}, \{\Sigma_i\})$ to $(f^*T_{\ov{\sX}} \To \sC \times S)$.  ($\Phi$ is defined on morphisms in the obvious way).  Since $\text{Pic}(\sC)$ is discrete \cite[Section 3.1]{Ca1}, $\Phi$ is locally constant. Thus the functions $(\sC \times S \xrightarrow{f} \ov{\sX}, \{\Sigma_i\}) \mapsto h^i(\sC, f^*T_{\ov{\sX}})$ are locally constant as well.  This proves the lemma.

We note that over the complex numbers the smoothness of $\sK(\sX)$ is asserted in \cite[Section 6.2]{AGV2}.
\end{proof}

\begin{no-num}{Lemma \ref{obs_theory}}
    \hfill
    \begin{enumerate}
        \item The natural map $$(R^{\bullet}\pi_* f^*T_{\sX})^{\vee} \ \xrightarrow{\phi} \ \Omega^1_{\sK(\sX)/ \sM_{0,3}^{tw}}$$ is a perfect relative obstruction theory with virtual dimension (denoted vdim) given by the locally constant function $$vdim  = dim \sX - \text{age}\circ e_1 - \text{age}\circ e_2 - \text{age}\circ e_3.$$
        \item $R^1\pi_* f^*T_{\sX}$ is locally free (denote the locally constant rank by $r$), and the virtual fundamental class (denoted $[\sK(\sX)]^{vir}$) in $A_{vdim}(\sK(\sX))_{\bQ}$ induced by $\phi$ is $$[\sK(\sX)]^{vir} = c_{r}(R^1\pi_*f^*T_{\sX}).$$
    \end{enumerate}
\end{no-num}
\begin{proof}
Part (2) follows from \cite[Prop. 7.3]{BF} using Part(1), Lemma \ref{K_smooth} above, and \cite[Thm. 5.2.1]{Kresch}.  (The fact that $R^1\pi_* f^*T_{\sX}$ is locally free also follows from the proof of Lemma \ref{K_smooth} above).  For Part (1), we proceed exactly as outlined in Section 4.5 of \cite{AGV2}.  Finally for the virtual dimension, it suffices to compute $\chi(\sC, f^*T_{\sX}) = \chi(\ov{\sC}, \ov{f}^*T_{\sX})$ where $(\sC \xrightarrow{f} \sX, \{\Sigma_i\})$ is a $\bF_q$-point in $\sK(\sX)$, and $(\ov{\sC} \xrightarrow{\ov{f}} \sX, \{\Sigma_i\})$ is the corresponding point in $\ov{\sK(\sX)}$.  But then the formula follows from the Riemann-Roch theorem on curves \cite[Thm. 7.2.1]{AGV2}.
\end{proof}

\begin{no-num}{Lemma \ref{K_compat}}
  \hfill
  \begin{enumerate}
    \item There is an equivalence
        $$\sK(\sX \times_{\bF_q} \ov{\bF}_q) \xrightarrow{\cong} \sK(\sX) \times_{\bF_q} \ov{\bF}_q.$$
    \item Under the identification above, the following functors are 2-isomorphic:
        $$\sK(F_{\ov{\sX}})\ \overset{\cong}{\rightrightarrows} \  F_{\ov{\sK(\sX)}}.$$
  \end{enumerate}
\end{no-num}
\begin{proof}
First we note that Part(1) appears in \cite[Prop. 5.2.1]{AV1}.  However, we include an elementary and explicit proof required for Part(2).

Denote the structure map by $A: \text{Spec}(\ov{\bF_q}) \To \text{Spec}(\bF_q)$, and denote projection maps by $p$ (e.g. $p_{\sX}: \ov{\sX} \To \sX$).  For part (1), an object of $\sK(\sX \times_{\bF_q} \ov{\bF}_q)$ over base $\ov{\bF}_q$-scheme $S$ is given by a commutative diagram
\begin{equation} \label{d1}
    \begin{diagram}
    \sqcup_{i=1}^3 \Sigma_i & \rInto & \sC  & \rTo^{f_S}            & \ov{\sX}          \\
                             &        & \dTo^{\pi_S} &                     &  \dTo        \\
                             &        & S    & \rTo^{q}      & \text{Spec}(\ov{\bF}_q)
    \end{diagram}
\end{equation}
where we have suppressed the map on coarse moduli spaces.  Define a functor $\phi$ by associating the pair consisting of the map $q$ together with the diagram given by composing $f_S$ (resp. $q$) in (\ref{d1}) with $p_{\sX}$ (resp. with $A$).  It is clear that $p_{\sX} \circ f_S$ is representable.  To see that it is stable, note that the map on coarse moduli schemes induced by $p_{\sX}$ is $p_{X}$, the projection onto $X$.  Moreover, $p_{X}^* \sO_{X}(1) \cong \sO_{\ov{X}}(1)$ and so $f_S$ is stable if and only if $p_{\sX} \circ f_S$ is stable.  One can define $\phi$ on morphisms in the obvious way, and it is easy to see that indeed $\phi$ is a functor.

For the reverse, an object of $\sK(\sX) \times_{\bF_q} \ov{\bF}_q$ over base $\ov{\bF}_q$-scheme $S$ is given by a commutative diagram
\begin{equation} \label{d2}
    \begin{diagram}
    \sqcup_{i=1}^3 \Sigma_i & \rInto & \sC  & \rTo^{f'_S}            & \sX          \\
                            &        & \dTo^{\pi_S} &               &  \dTo        \\
                            &        & S    & \rTo^{q'}              & \text{Spec}(\bF_q)
    \end{diagram}
\end{equation}
where $q'$ factors through the structure map $q'': S \To \text{Spec}(\ov{\bF}_q)$ (i.e. $A \circ q'' = q'$).  Then define $\psi$ (inverse of $\phi$) by sending (\ref{d2}) to the diagram given by (\ref{d1}) with $q$ replaced by $q''$ and $f_S$ replaced by the unique map induced by the pair $(q'' \circ \pi_S, f'_S)$.  One can define $\phi$ on morphisms in the obvious way.  It easy to check that indeed $\psi$ is a functor and that the pair $(\phi, \psi)$ gives the required equivalence.

For the second part, we show that the following diagram is 2-commutative.
\begin{diagram}
\ov{\sK(\sX)} & \rTo^{F_{\ov{\sK(\sX)}}} & \ov{\sK(\sX)} \\
\uTo_{\phi} & & \dTo_{\psi} \\
\sK(\ov{\sX}) & \rTo^{\sK(F_{\ov{\sX}})} & \sK(\ov{\sX})
\end{diagram}
Let $\eta$ denote the diagram (\ref{d1}) in $\sK(\ov{\sX})$ above.  Since $F_{\ov{\sK(\sX)}} = 1_{\sX} \times F_{\text{Spec}(\ov{\bF}_q)}$, we have that $\phi \circ F_{\ov{\sK(\sX)}} \circ \psi (\eta)$ is given by a diagram similar to (\ref{d2}) with $q$ replaced by $F_{\text{Spec}(\ov{\bF}_q)} \circ q$ and $f_S$ replaced by the unique morphism induced by the pair $(F_{\text{Spec}(\ov{\bF}_q)} \circ q \circ \pi_S, p_{\sX} \circ f_S)$.  Denote this unique morphism by $\ast$.  On the other hand, consider the commutative diagram
\begin{diagram}
\sC_S       & \rTo^{f_S}    & \ov{\sX}      & \rTo^{F_{\ov{\sX}}}   & \ov{\sX}      & \rTo^{p_{\sX}}    & \sX \\
\dTo^{\pi_S}&               & \dTo          &                       & \dTo          &                   & \dTo \\
S           & \rTo^{q}      & \text{Spec}(\ov{\bF}_q)    & \rTo^{F_{\text{Spec}(\ov{\bF}_q)}} & \text{Spec}(\ov{\bF}_q) & \rTo^{A}  & \text{Spec}(\bF_q)
\end{diagram}
with two right squares Cartesian and $F_{\ov{\sX}} \circ p_{\sX} = p_{\sX}$.  This diagram shows that $\ast$ is given by $F_{\ov{\sX}} \circ f_S$.  Thus we see that $\phi \circ F_{\ov{\sK(\sX)}} \circ \psi (\eta) = \sK(F_{\ov{\sX}})(\eta)$.  One can further check that $\phi \circ F_{\ov{\sK(\sX)}} \circ \psi$ and $\sK(F_{\ov{\sX}})$ agree on morphism as well and this proves the lemma.
\end{proof}

\begin{no-num}{Lemma \ref{F_on_vir_class}}
The equation $$\sK(F_{\ov{\sX}})_* [\sK(\ov{\sX})]^{vir} = q^{-vdim}[\sK(\ov{\sX})]^{vir}$$ holds in $H_*(\sK(\ov{\sX}), \bQ_l)$.
\end{no-num}
\begin{proof}
Consider the Cartesian diagram
\begin{diagram}
\sK(\ov{\sX})   & \rTo^{p_{\sK(\sX)}}   & \sK(\sX) \\
\dTo            &                       &  \dTo \\
\text{Spec}(\ov{\bF}_q) & \rTo^{A}      & \text{Spec}(\bF_q) \\
\end{diagram}
where we identify $\ov{\sK(\sX)}$ with $\sK(\ov{\sX})$ by Lemma \ref{K_compat}.  The projection $p_{\sK(\sX)}$ induces a map on Chow groups $p_{\sK(\sX)}^*: A_*(\sK(\sX))_{\bQ} \To A_*(\sK(\ov{\sX}))_{\bQ}$.  We first show that $p_{\sK(\sX)}^* [\sK(\sX)]^{\text{vir}} = [\sK(\ov{\sX})]^{\text{vir}}$. Let $\sU_{\sX}$ (resp. $\sU_{\ov{\sX}}$) denote the universal curves over $\sK(\sX)$ (resp. $\sK(\ov{\sX})$).  By \cite[Corollary 9.1.3]{AV1}, $\sU_{\sX}$ (resp. $\sU_{\ov{\sX}}$) is an open and closed substack of $\sK_{0,4}(\sX, 0)$ (resp. $\sK_{0,4}(\ov{\sX}, 0))$.  Thus the proof of Lemma \ref{K_compat} identifies $\sU_{\ov{\sX}}$ with $\sU_{\sX} \times_{\bF_q} \ov{\bF}_q$, and we have the following diagram with Cartesian squares:
\begin{diagram}
                       &        & \ov{\sX}                  & \rTo^{p_{\sX}}    &                   &\sX  \\
                       & \ruTo^{f_{\ov{\sX}}} &             &                   & \ruTo^{f_{\sX}}   & \\
\sU_{\ov{\sX}}         &        & \rTo^{p_{\sU_{\sX}}}      & \sU_{\sX}         &                   &\\
\dTo^{\pi_{\ov{\sX}}}  &        &                           & \dTo_{\pi_{\sX}}  &             &     \\
\sK(\ov{\sX})          &        & \rTo^{p_{\sK(\sX)}}       & \sK(\sX).          &             &     \\
\end{diagram}
Since $p_{\sX}^*T_{\sX} \cong T_{\ov{\sX}}$ and $p_{\sK(\sX)}^* R^1(\pi_{\sX})_* \cong R^1(\pi_{\ov{\sX}})_* p_{\sU_{\sX}}^*$ a simple diagram chase gives
\begin{alignat}{1}
p_{\sK(\sX)}^* [\sK(\sX)]^{\text{vir}}  &= c_{vdim}(R^1(\pi_{\sX})_*f_{\sX}^*T_{\sX}) \notag \\
                                        &= c_{vdim}(R^1(\pi_{\ov{\sX}})_*f_{\ov{\sX}}^*T_{\ov{\sX}}) \notag \\
                                        &= [\sK(\ov{\sX})]^{\text{vir}}. \notag
\end{alignat}

Now fix a connected component $\sK(\sX)_0 \subset \sK(\sX)$ so that the virtual dimension (vdim) is constant.  If the virtual class on $\sK(\sX)_0$ is given by $\sum n_i [V_i]$ where each $V_i$ is a substack of pure dimension vdim, then $[\sK(\ov{\sX})]^{\text{vir}} = \sum n_i [\ov{V}_i]$.  Thus we have $\sK(F_{\ov{\sX}})(\ov{V}_i) = F_{\sK(\ov{\sX})}(\ov{V}_i) = \ov{V}_i$.  Furthermore if $F_{\sK(\ov{\sX})}^{\text{geo}}$ denotes the geometric Frobenius morphism, then we also have $F_{\sK(\ov{\sX})}^{\text{geo}}(\ov{V}_i) = \ov{V}_i$.  Hence $(F_{\sK(\ov{\sX})}^{\text{geo}})_*\ [\ov{V}_i] = q^{vdim} [\ov{V}_i]$ in the Chow group.  But after passing to (co)homology the arithmetic and geometric Frobenius maps are inverses,
and this proves the lemma.
\end{proof}

\bibliographystyle{hamsplain}
\bibliography{xbib}
\end{document}